\documentclass [a4paper,twoside,12pt]{article}

\usepackage[utf8]{inputenc}
\usepackage{amsmath,amsfonts,amssymb}
\usepackage{vmargin,graphicx,theorem}
\usepackage[english]{babel}
\usepackage{enumerate}
\usepackage{color}
\usepackage{pstricks}
\usepackage{mdwlist}
\usepackage{pst-fill,pst-grad,pst-plot,pst-eucl,pstricks-add,pst-node}

\setpapersize[portrait]{A4}
\setmarginsrb{1.5cm}{1cm}{1.5cm}{2.5cm}{1.5cm}{0cm}{0.5cm}{2cm}

\newcommand{\rr}{\ensuremath{\mathbb{R}}}

\newcommand{\dE}{\ensuremath{\mathbb{E}}}

\newcommand{\dR}{\ensuremath{\mathbb{R}}}
\newcommand{\dS}{\ensuremath{\mathbb{S}}}

\newtheorem{ethm}{Theorem}[section]

\newtheorem{ecor}[ethm]{Corollary}

\newtheorem{eprop}[ethm]{Proposition}

\newtheorem{elem}[ethm]{Lemma}

\newtheorem{edefi}[ethm]{Definition}

\newtheorem{erem}[ethm]{Remark}

\newcommand{\proofend}{~$\rhd$}
\newcommand{\proofbegin}{~$\lhd$}

\newenvironment{eproof}
               {\noindent {\emph{\textbf{Proof}}}\\\proofbegin~}
               {\proofend\\}

\newcommand{\ABS}[1]{\ensuremath{{\left| #1 \right|}}} 
\newcommand{\PAR}[1]{\ensuremath{{\left(#1\right)}}} 
\newcommand{\BRA}[1]{\ensuremath{{\left\{#1\right\}}}} 

\renewcommand{\phi}{\varphi}

\renewcommand{\geq}{\geqslant}

\newcommand{\R}{\dR}

\newcommand{\n}{\nabla}

\newcommand{\beq}{\begin{equation}}\newcommand{\eeq}{\end{equation}}
\begin{document}

\title{ Convergence to equilibrium in Wasserstein distance for Fokker-Planck equations}

\author{ Fran{\c c}ois Bolley\thanks{Ceremade, Umr Cnrs 7534, Universit\'e Paris-Dauphine, Place de Lattre de Tassigny, F-75016 Paris cedex 16. bolley@ceremade.dauphine.fr}, Ivan Gentil\thanks{Institut Camille Jordan, Umr Cnrs 5208, Universit\'e Claude Bernard Lyon 1, 43 boulevard du 11 novembre 1918,
F-69622 Villeurbanne cedex. gentil@math.univ-lyon1.fr}\, and Arnaud Guillin\thanks{Institut Universitaire de France and Laboratoire de Math\'ematiques, Umr Cnrs 6620, Universit\'e Blaise Pascal, Avenue des Landais, F-63177 Aubi\`ere cedex. guillin@math.univ-bpclermont.fr}}

\date{July 6, 2012}

\maketitle

\abstract{We describe conditions on non-gradient drift diffusion Fokker-Planck equations for its solutions to converge to equilibrium with a uniform exponential rate in Wasserstein distance. This asymptotic behaviour is related to a functional inequality, which links the distance with its dissipation and ensures a spectral gap in Wasserstein distance. We give practical criteria for this inequality and compare it to  classical ones.  The key point is to quantify the contribution of the diffusion term to the rate of convergence, in any dimension, which to our knowledge is~a~novelty.}
\bigskip

\bigskip

 \section*{Introduction}

In this work we consider the Fokker-Planck equation
 \begin{equation}\label{eq:FP}
 {\partial_t \mu_t} = \nabla \cdot (\nabla \mu_t + \mu_t A), \qquad t>0, x \in \rr^n
 \end{equation}
where $A$ is given vector field on $\rr^n$. The evolution preserves mass and positivity, and we are concerned with initial data $\mu_0$ which are probability measures on $\rr^n$, so that so are the solutions $\mu_t = \mu(t,.)$ at any time $t>0$.  For measures on $\rr^n$ with a density with respect to the Lebesgue measure, we shall use the same notation for the measure and its density, hoping that it is not confusing.

We are interested in criteria ensuring uniform bounds on the long time behaviour of solutions.
 
 \medskip
To explain our main issue, let us start with the classical case when $A=\nabla V$ with $\int e^{-V} dx =~1$.  
 The probability measure $d\nu(x) = e^{-V(x)} dx$ is a stationary solution of~\eqref{eq:FP} and it is interesting to know for which $V$ all solutions $\mu_t$ converge to $\nu$ as $t$ tends to infinity, in which sense and with a rate.

 There are various ways of measuring the gap between a solution $\mu_t$ of the equation and the stationary one $e^{-V}$: total variation (as in Meyn-Tweedie's approach), $L^2$-norm, relative entropy, Wasserstein distance. Perhaps, the simplest way is to consider the  $L^2$-norm 
 $$
 G(t)=\int \PAR{\mu_t-e^{-V}}^2 \, e^{V}dx = \int \PAR{\frac{\mu_t}{e^{-V}}-1 }^2 \, e^{-V}dx
 $$
of the difference of the densities, with the weight $e^V$. Formally, by integration by parts,
$$
G'(t)
=2\int \Big( \frac{\mu_t}{e^{-V}}-1 \Big) {\partial_t \mu_t} dx
=-2\int \ABS{\nabla\PAR{\frac{\mu_t}{e^{-V}} -1} }^2 \, e^{-V} \, dx, \qquad t>0.
$$
Here $\vert \cdot \vert$ is the Euclidean norm on $\rr^n$. In particular the quantity $G(t)$ is non-increasing in time. 

Assume now that the measure $e^{-V}$ satisfies a  Poincar\'e inequality with constant $C>0$, that~is, 
\begin{equation}\label{eq-poin}
\int \PAR{f-\int f e^{-V}}^2 e^{-V} \, dx\leq \frac{1}{C} \int |\nabla f|^2 e^{-V} \, dx
\end{equation}
for all $f$. By choosing $f = \mu_t / e^{-V}-1$,  we obtain  $G'(t)\leq - 2 C G(t)$. Hence
\begin{equation}\label{eq:contL21sol}
\int|\mu_t-e^{-V}|^2 \, e^{V}dx\leq e^{-2 C t}\int|\mu_0-e^{-V}|^2 \, e^{V}dx, \qquad t \geq 0
\end{equation}
by integration. In particular this ensures the strong convergence of $\mu_t$ to $e^{-V}$ in $L^2(e^V)$ for any initial datum $\mu_0$ in $L^2(e^V)$. In fact, \eqref{eq:contL21sol} is equivalent to~\eqref{eq-poin} by time-differentiating at $t=0$. 

Then simple criteria are known for a measure $e^{-V}$ to satisfy the Poincar\'e inequality~\eqref{eq-poin}: for instance,~\eqref{eq-poin} holds  if the Hessian matrix $\nabla^2 V(x)$ is uniformly bounded by below by $C {\rm Id}_n$ (Bakry-\'Emery criterion, see~\cite{logsob} for instance); more generally it holds for a $C$ if $V$ is convex, see for example \cite{bbcg08}. The argument can also be performed for diverse convex functionals of the quantity $\mu_t/e^{-V}$, under the name of entropy method (see~\cite{amtucpde01} for instance).

In fact the Poincar\'e inequality~\eqref{eq-poin} implies the following stronger contraction property between any two solutions: if $\mu_t$ and $\nu_t$ are two  solutions in $L^2(e^V)$, then~\eqref{eq-poin} with $f = (\mu_t - \nu_t)/e^{-V}$ leads to
\begin{equation}\label{eq:contL22sol}
\int|\mu_t-\nu_t|^2 \, e^{V} dx \leq e^{-2Ct}\int|\mu_0-\nu_0|^2 \, e^{V} dx, \qquad t \geq 0.
\end{equation}
 It implies~\eqref{eq:contL21sol} by letting $\nu_0 = e^{-V}$.

 As a conclusion, the long time convergence estimate~\eqref{eq:contL21sol} is equivalent to the (seemingly stronger) $L^2$-contraction property~\eqref{eq:contL22sol} of two solutions,  and to the Poincar\'e inequality~\eqref{eq-poin}. 
 
 \bigskip
 
 Contraction results between solutions to~\eqref{eq:FP} can also be measured in terms of Wassertein distances. If $\rho_1$, $\rho_2$ are two probability measures on $\rr^n$, their Wasserstein distance is defined~by 
$$
W_2(\rho_1,\rho_2)=\inf \big( \dE \vert X-Y \vert^2 \big)^{1/2}, 
$$
where the infimum runs over all random variables $X$ and $Y$ with law respectively $\rho_1$ and $\rho_2$. This distance metrizes a weak convergence (as opposed to the strong $L^2$ convergence above), but has the advantage of being defined on the larger and more natural space of probability measures on~$\rr^n$. Moreover convergence for this distance can be turned into convergence in Sobolev norms by means of interpolation estimates as in~\cite{portoercole}. It is adapted to~\eqref{eq:FP} since, by the It\^o formula, a measure solution $\mu_t$  to~\eqref{eq:FP} can be seen as the law at time $t$ of the process $(X_t)_{t \geq 0}$ solution to  the stochastic differential equation 
 \begin{equation}
 \label{eq-sde}
 dX_t=\sqrt{2}\,dB_t-\nabla V(X_t)dt.
 \end{equation}
 Here $(B_t)_{t\geq 0}$ is a standard Brownian motion in $\dR^n$ and the initial datum $X_0$ has law $\mu_0$. 
  
Let now $\mu_0$ and $\nu_0$ be two measures on $\R^n$, and  $(X_t)_{t}$ (resp. $(Y_t)_{t}$) the solution to~\eqref{eq-sde} starting from $X_0$ of law $\mu_0$ (resp. $Y_0$ of law $\nu_0$), both driven by the same Brownian motion. Then 
 $$
 \frac{d}{dt} |X_t-Y_t|^2= - 2 \, \left(\nabla V(X_t)-\nabla V(Y_t)\right) \cdot (X_t - Y_t).
 $$
 Now, if $V$ satisfies  $\nabla^2 V (x)\geq C {\rm Id}_n$ for all $x\in\dR^n$ and for a $C\in\dR$, that is,
 \begin{equation}
\label{eq-hypo-a}
 (\nabla V(x) - \nabla V(y)) \cdot  (x-y) \geq C \,  \vert x-y \vert^2 
\end{equation}
for all $x,y\in\R^d$, then 
 $$
\dE \, \vert X_t-Y_t \vert^2 \leq e^{-2C t} \, \dE \, \vert X_0-Y_0 \vert^2
 $$
by integrating in time and taking the expectation. Moreover
$$
W_2^2(\mu_t,\nu_t)\leq  \dE \vert X_t-Y_t \vert^2
$$
since $X_t$ and $Y_t$ have respective laws $\mu_t$ and $\nu_t$.  Now taking the infimum over  $X_0$ and $Y_0$ gives the following contraction-type estimate between any two solutions  : 
\begin{equation}
\label{eq-contraction-wasserstein}
W_2(\mu_t,\nu_t)\leq e^{- Ct}W_2(\mu_0,\nu_0), \qquad t \geq 0.
\end{equation}
Such a contraction-type estimate is a key estimate in the theory of gradient flows in the space of probability measures, an instance of which is~\eqref{eq:FP} when $A = \nabla V$ (see~\cite{ambrosio-gigli-savare}).

In particular, by choosing $\nu_0$ as the stationary solution $e^{-V}$ it implies the bound 
\begin{equation}
\label{eq-wasserstein-decay}
W_2(\mu_t, e^{-V})\leq e^{-Ct}W_2(\mu_0, e^{-V}), \qquad t \geq 0
\end{equation}
 for any initial condition $\mu_0$.  For $C>0$ it ensures that $e^{-V}$ is the only stationary state of~\eqref{eq:FP} and quantifies the convergence of all solutions to it; it can be seen as a spectral gap in Wasserstein distance. 
 
 \medskip
 
 Of course~\eqref{eq-contraction-wasserstein} is a stronger statement than~\eqref{eq-wasserstein-decay} since it enables to compare any two solutions, and not only a solution to the stationary one. But it asks for extremely strong assumptions on the drift: indeed, according to K.-T. Sturm and M. von Renesse, the uniform convexity condition~\eqref{eq-hypo-a}  is in fact {\it equivalent} to~\eqref{eq-contraction-wasserstein}; more generally  when the vector field $A$ is not necessary a gradient, then solutions of~\eqref{eq:FP} satisfy~\eqref{eq-contraction-wasserstein}   if and only if~\eqref{eq-hypo-a} holds with $A$ instead of $\nabla V$ (see \cite{sturm-vonrenesse} and 
Remark~\ref{rem:sturmvonrenesse}, 
 and also~\cite{nps} for a duality proof of the sufficient condition). In this case, and with $C>0$, this classically ensures the existence and uniqueness of a stationary solution, as used in diverse contexts in~\cite{bgm},~\cite{portoercole}, or~\cite{cmcv-06} for instance.

\bigskip

The purpose of this work is twofold: First, to consider possibly non-gradient drifts $A$, which naturally appear for example in polymeric fluid flow or Wigner-Fokker-Planck equation (see \cite{jllo06} or \cite{acm10}).  Such non gradient drifts forbid the gradient flow approach to~\eqref{eq:FP}, which holds only in the gradient case. Then, and above all, to give weaker conditions than~\eqref{eq-hypo-a} on the drift $A$ for the uniform convergence estimate~\eqref{eq-wasserstein-decay} to hold for solutions to~\eqref{eq:FP}.  As for the $L^2$-norm and the Poincar\'e inequality, it will be described by a functional inequality, which links the Wasserstein distance with its dissipation along the flow of the equation. As will be seen later on, an interesting fact is that it holds for potentials which are uniformly convex only at infinity. For that purpose we will use the diffusion term to overcome the possible degeneracy of the potential convexity in some region. We will see on examples how an a priori polynomial rate of convergence can simply be turned into an exponential rate by this method. To our knowledge this is the first quantitative use of the contribution of the diffusion term in measuring the convergence to equilibrium in Wasserstein distance in any dimension (this idea also appears in~\cite{carrillodft} in the $1d$ case, and with a crucial use of the specific $1d$ formulation of the distance).\\

In Section~\ref{sec-frame}, we introduce the objects studied in the paper. In Section~\ref{sec: convergence} we derive the Wasserstein distance dissipation along solutions to~\eqref{eq:FP} when $A$ is not necessarily a gradient, and state first simple criteria for the uniform stability or convergence estimate~\eqref{eq-wasserstein-decay}. In Section~\ref{sec-wj} we introduce the $WJ$ inequality which governs~\eqref{eq-wasserstein-decay}, and give further practical conditions to this inequality and its connections with classical functional inequalities as the Poincar\'e or logarithmic Sobolev inequalities. 

\bigskip

 Let us finish by some possible extension to nonlinear models. For example, contraction properties such as~\eqref{eq-contraction-wasserstein} also hold for nonlinear equations such as the granular media equation
$$
\partial_t \mu_t = \nabla \cdot (\nabla \mu_t + \mu_t (\nabla V + \nabla W \ast_x \mu_t)), \qquad t>0, x \in \rr^n
$$ 
under hypothesis like~\eqref{eq-hypo-a} on the potentials $V$ and $W$ (see~\cite{cmcv-06}); here $\ast_x$ stands for the convolution on $\rr^n$. It is then natural to hope that we can go beyond this strict convexity assumption using our approach. In~\cite{bgg} we precisely show that the method is sufficiently robust to include non-uniformly convex potentials.

\section{Framework}
\label{sec-frame}

We consider the Fokker-Planck equation starting from a probability measure $\mu_0$, \begin{equation}
\label{eq-fp1}
\partial_t \mu_t =\nabla\cdot(\nabla \mu_t+\mu_t A)=\nabla\cdot(\mu_t ( \nabla\log \mu_t+A)), \qquad t>0, x \in \rr^n
\end{equation}
where $A$ is a $\mathcal C^1$ function on $\dR^n$ and $\nabla\cdot G$ is  the divergence of a vector field $G$.

The existence of a non-explosive solution can be proven under simple conditions on $A$. For instance, if there exist $a$ and $b$ such that
$$
x\cdot A (x)\geq -a|x|^2-b
$$
for all $x$, then for any initial datum $\mu_0$ in the space $\mathcal P_2(\rr^n)$ of probability measures $\rho$ on $\rr^n$ such that $\int |x|^2d\rho(x) <\infty$ there exists a continuous curve $(\mu_t)_{t \geq 0}$ of probability measures such that~\eqref{eq-fp1} holds in the sense of distributions. We shall assume that for any $t>0$ a solution $\mu_t$ is in $\mathcal P_2(\rr^n)$ and has a $\mathcal C^1$ positive density with respect to the Lebesgue measure: this is proven in diverse frameworks for instance in~\cite{stroock}, the appendix in~\cite{bgv07}, Corollary 3.6 in~\cite{bogachev08}, see also~\cite{nps} and the references therein.

It\^o's formula  implies that the law $(\mu_t)_{t\geq0}$ 
of the Markov process 
\begin{equation}\label{eq-eds}
dX_t=\sqrt{2}dB_t- A(X_t)dt,
\end{equation}
where $X_0$ has law $\mu_0$ and $(B_t)_{t\geq0}$ is a Brownian motion on $\dR^n$,  is a solution to~\eqref{eq-fp1}.  Equation~\eqref{eq-fp1} is also called the Kolmogorov forward equation.

We assume that there exists a positive smooth stationary solution  $e^{-V}$ of~\eqref{eq-fp1}, which is a probability measure and where $V$ is a $C^2$ map on $\rr^n$. Letting $F = A - \nabla V$, equation~\eqref{eq-fp1} reads
\begin{equation}
\label{eq-fp}
\partial_t \mu_t =\nabla\cdot(\mu_t(\nabla\log \mu_t+\nabla V+F)).
\end{equation}
Here the vector field $F$ satisfies $\nabla\cdot(e^{-V}F)=0$, which is a necessary and sufficient condition for $e^{-V}$ to be a stationary solution.

Let 
$\nabla A$ be the Jacobian matrix of $A$ and $\nabla^S A = ( 
\nabla A + \nabla A 
^T )/2$ be its symmetric part. 
We saw in the introduction that the condition $\nabla^S A  \geq C \, {\rm Id}_n$ as symmetric matrices on $\rr^n$, 
with $C>0$ and uniformly on $\rr^n$, ensures the existence of a unique stationary solution in the space of probability measures, and convergence of all solutions to it. Weaker conditions on $A$ for such an existence can be obtained by Liapunov methods for instance, but deriving quantitative estimates on the steady state, and a fortiori convergence estimates, only from the knowledge of $A$, is an interesting and difficult issue, which will not be addressed in the present work. We refer to~\cite{arnoldcarlen} for an entropy dissipation approach to convergence rates, and with a general diffusion matrix.

\medskip
 
The generator  $L^*$ defined by $L^*f=\Delta f +\nabla\cdot(f(\nabla V+F))$ for $f$ a $C^2$ map on $\rr^n$ is the dual operator in $L^2(dx)$ of $L$ defined by $Lf=\Delta f-\nabla f\cdot(\nabla V+F)$. Moreover  $L$ is the infinitesimal generator of the Markov semigroup $(P_t)_{t\geq0}$ defined by
$$
P_t f(x)=E_x(f(X_t))
$$
for any smooth function $f$; here $(X_t)_{t\geq0}$ is the Markov process, solution of the stochastic differential equation~\eqref{eq-eds}, such that $X_0=x$. In other words, the function $P_tf$ solves the partial differential equation 
\begin{equation}
\label{eq-ou}
\partial_t u=Lu,
\end{equation}
with initial datum $f$. 

\medskip

If $\mu_t$ is a solution to~\eqref{eq-fp} then $\varphi_t=e^{V} \mu_t$  satisfies the PDE
\begin{equation}
\label{eq-ouT}
\partial_t \varphi_t=\Delta \varphi_t-\nabla\varphi_t \cdot(\nabla V-F).
\end{equation}
Conversely, if $\varphi_t$ is a smooth positive solution to~\eqref{eq-ouT} with initial datum $\varphi_0$ such that $\int \varphi_0 e^{-V} dx=1$, then 
$$
\mu_t = e^{-V} \varphi_t
$$
for $t \geq 0$ is a positive probability density which solves~\eqref{eq-fp} with the initial datum $\varphi_0 e^{-V}$. 
The diffusion operator $L^\top f =\Delta f-\nabla f\cdot(\nabla V-F)$ can now be seen as the infinitesimal  generator of a Markov semigroup denoted $(P^\top_t)_{t\geq 0}$. It is the dual of $L$ in $L^2(d\nu)$, where $d\nu=e^{-V}dx$, that~is, 
$$
\int fLg\,d\nu=\int gL^\top f\,d\nu
$$  
for all compactly supported $C^2$ functions $f$ and $g$.

\medskip

Moreover, the measure  $d\nu=e^{-V}dx$ is an invariant measure for both generators $L$ and $L^\top$, that is, for all compactly supported $C^2$ functions $f$
$$
\int L^\top f d\nu=\int L f d\nu=0.
$$

\medskip

When $A=\nabla V$ (or equivalently $F=0$), then~\eqref{eq-fp} is the usual Fokker-Planck equation whereas the dual form~\eqref{eq-ou} is the general Ornstein-Uhlenbeck equation. In that case $L^\top=L$ and $L$ is symmetric in $L^2(d\nu)$ : we say that $\nu$ is reversible. 
\medskip

The discrepancy between probability measures will mainly be estimated in terms of the Wasserstein distance:
or two measures $\rho_1$ and $\rho_2$ in $\mathcal P_2(\rr^n)$ it is defined by 
$$
W_2(\rho_1,\rho_2)=\inf \Big( \int_{\rr^{2n}}  |x-y|^2d\pi(x,y) \Big)^{1/2},
$$
where the infimum runs over all probability measures $\pi$ on $\dR^n\times\dR^n$ with marges $\rho_1$ and $\rho_2$, that is, for all bounded functions $f$ and $g$ on $\rr^n$
$$
\int_{\rr^{2n}} (f(x)+g(y)) \, d\pi(x,y)=\int_{\rr^{n}} fd\rho_1+\int_{\rr^{n}} gd\rho_2
$$ 
(see~\cite{ambrosio-gigli-savare} or~\cite{villani-book1} for example). This definition is of course the same as the one given in the introduction in terms of random variables. All the measures considered in the sequel will be in $\mathcal P_2(\rr^n)$, even if not specified.

\medskip

Brenier's Theorem gives an explicit expression of the Wasserstein distance: if $\rho_1$ is absolutely continuous with respect to the Lebesgue measure then there exists a convex function $\phi$ such that $\nabla \phi\#\rho_1=\rho_2$, that is, 
$$
\int_{\rr^{n}} g \, d\rho_2=\int_{\rr^{n}} g(\nabla \phi) \,d\rho_1
$$
for every bounded test  function $g$; moreover
$$
W_2^2(\rho_1,\rho_2)=\int_{\rr^{n}} |x-\nabla\phi(x)|^2 \, d\rho_1(x).
$$
\medskip

The Legendre transform will be useful for the next sections: for a map $\phi:\dR^n\mapsto\dR\cup\{\infty\}$ it is the map $\phi^*:\dR^n\mapsto\dR\cup\{\infty\}$ defined by 
$$
\phi^*(q)=\sup_{x\in\dR^n}\BRA{q\cdot x-\phi(x)}. 
$$
If $\rho_1$ and $\rho_2$ are probability densities in $\mathcal P_2(\rr^n)$ such that $\n \phi \#\rho_1 = \rho_2$, then $\n \phi^* \#\rho_2 = \rho_1$.

\section{Convergence in Wasserstein distance}
\label{sec: convergence}

Convergence in Wasserstein distance is related to its time-derivative, which was studied by L.~Ambrosio, N.~Gigli and G.~Savar\'e in~\cite[Th. 8.4.7]{ambrosio-gigli-savare} (see also~\cite[Th. 23.9]{villani-book1}). 

For a probability measure $\rho_1$ and a probability density $h$ with respect to $\rho_1$ we let
\begin{equation}
\label{eq-fisher}
H(\rho_2 \vert \rho_1) = \int h \, \log h \, d\rho_1, \qquad 
I(\rho_2 \vert \rho_1)
=
\int \frac{\vert \nabla h \vert^2}{h}  \, d\rho_1
\end{equation}
respectively be the entropy and the Fisher information  of $\rho_2 = h \rho_1$ with respect to $\rho_1$.

\begin{ethm}[\cite{ambrosio-gigli-savare}]
\label{thm-derivation}
Assume that $V,F$ 

are such that $\int |F|^4d\nu<\infty$ with $d\nu = e^{-V} dx \in \mathcal P_2(\rr^n)$. 
Let  $\mu_t$ be a solution of~\eqref{eq-fp} with initial condition having a smooth density $\mu_0 \in \mathcal P_2(\rr^n)$ such that  
\begin{equation}
\label{eq-regularity}
\int \mu_0^2 e^{V}dx<\infty.
\end{equation}

Then the map $t\mapsto W_2(\mu_t,\nu)$ is absolutely continuous and for  almost every $t \geq 0$
\begin{equation}
\label{eq-derivation}
\frac12 \frac{d}{d t}W_2^2(\mu_t,\nu)=\int \PAR{\nabla \psi_t-x}\cdot \PAR{\nabla \log\mu_t+A}d\mu_t
\end{equation}
where for every $t\geq 0$, ${\nabla \psi_t}\#{\mu_t}=\nu$.  
\end{ethm}

Let us first give a direct and formal proof of this result.  Brenier's Theorem implies that  
$$
W_2^2(\mu_t,\nu)=\int \ABS{\nabla \phi_t(x) - x}^2d\nu(x)
$$
for all $t\geq0$, where $\nabla \phi_t \# \nu = \mu_t$. Then by formal time-differentiation
$$
\frac12 \frac{d}{d t}W_2^2(\mu_t,\nu)=\int \PAR{\nabla \phi_t(x) - x}\cdot\partial_t\nabla\phi_t d\nu.
$$ 
Now for all $g=g(t,x)$ the time-derivative of $\int g(\nabla\phi_t)d\nu=\int gd{\mu_t}$ is
$$
\int \nabla g(t, \nabla\phi_t) \cdot \partial_t \nabla \phi_t \, d\nu=\int g \, d(\partial_t\mu_t).
$$
For $g(t, x)=\frac{|x|^2}{2} - \phi^*_t (x)$, which satisfies $\nabla g(t, \nabla\phi_t (x))=\nabla\phi_t (x) -x$ by Legendre transform properties, this gives 
$$
\frac12 \frac{d}{d t}W_2^2(\mu_t,\nu) =\int \PAR{\frac{|x|^2}{2} - \phi_t^*} d(\partial_t\mu_t).
$$
 An integration by parts implies~\eqref{eq-derivation} with $\psi_t = \phi_t^*$. 
 
 \medskip
 
Another approach, developed in~\cite[Th. 4.1]{AGS11}, goes as follows : for given $t$ let $g(x) = \frac{\vert x \vert^2}{2} - \varphi_t^*(x)$ as above and $\bar{g}(y) = \frac{\vert y \vert^2}{2} - \varphi_t (y)$ be the Kantorovich potentials such that $g(x) + \bar{g}(y) \leq \frac{1}{2} \vert x-y \vert^2$ for every $x,y$ and 
$$
\frac{1}{2} W_2^2(\mu_t, \nu) = \int g \, d\mu_t + \int \bar{g} \, d\nu.
$$
We observe that 
$$
\frac{1}{2} W_2^2(\mu_{t-h}, \nu) \geq \int g \, d\mu_{t-h} + \int \bar{g} \, d\nu
$$
so that taking the difference, dividing by $h>0$ and letting $h \to 0^+$
$$
 \frac{1}{2} \frac{d}{dt} W_2^2(\mu_t, \nu) \leq \int g \, d (\partial_t \mu_t)
$$
as above.

\medskip

{\medskip

\begin{eproof}
It is a direct application of~\cite[Th.~23.9]{villani-book1} and we now check its assumptions.

First, the  vector field $\xi_t=\nabla\log\mu_t+\nabla V+F$ is locally Lipschitz since the solution $\mu_t$ has a smooth and positive density on $(0,\infty)$. Let us now check that 
$$
\int_{t_1}^{t_2}\int |\xi_t|^2d\mu_t \, dt<\infty
$$ 
for every $0<t_1<t_2$. Indeed
$$
\int |\xi_t|^2d\mu_t=\int \Big|\nabla\log \frac{\mu_t}{\nu} +F \Big|^2d\mu_t \leq  2 \, I(\mu_t  \vert \nu) + 2\int |F|^2d\mu_t. 
$$

On the one hand, since $\nabla \cdot (F e^{-V}) =0$,
$$
\int_{t_1}^{t_2} I(\mu_t \vert \nu) \, dt \leq \int_{0}^{t_2} I(\mu_t \vert \nu) \, dt = H(\mu_0 \vert \nu) - H(\mu_{t_2} \vert \nu) \leq H(\mu_0 \vert \nu)
$$
which is finite since so is $\int \mu_0^2 e^V dx .$

As for the other term, by the Cauchy-Schwarz inequality, 
\begin{multline*}
\int |F|^2d\mu_t=\int |F|^2 \, \frac{\mu_t}{\nu} \, d\nu \leq \PAR{\int |F|^4d\nu}^{1/2} \PAR{\int \left(\frac{\mu_t}{\nu}\right)^2 \, d\nu}^{1/2}
\\
\leq \PAR{\int |F|^4d\nu}^{1/2}\PAR{\int  \left(\frac{\mu_0}{\nu}\right)^2 \, d\nu}^{1/2}.
\end{multline*}
The last two bounds imply 
$$
\int_{t_1}^{t_2}\int |\xi_t|^2d\mu_tdt\leq 2\, H(\mu_0 \vert \nu)+ 2 (t_2-t_1)\PAR{\int |F|^4d\nu \int \mu_0^2e^{V}dx}^{1/2}<\infty.
$$ 
 \end{eproof}

\begin{erem}
The assumptions of Theorem~\ref{thm-derivation} hold for instance for the couple $(F, V)$ considered in~\cite{arnoldcarlen}.
\end{erem}

\begin{erem}
In the (gradient flow) case when $F=0$, the proof above requires the weaker condition $H(\mu_0 \vert \nu) < \infty$ instead of~\eqref{eq-regularity}. In fact, it is observed in the Theorem in~\cite{ov01} that $H(\mu_{t_1} \vert \nu) < \infty$ for any $t_1 >0$ and initial datum $\mu_0 \in \mathcal P_2(\rr^n)$ if $\nabla^2 V $ is uniformly bounded from below (by a possibly negative constant $\lambda$), so that Theorem~\ref{thm-derivation} can be extended to all solutions with initial datum in $\mathcal P_2(\rr^n)$; this is a general feature of gradients flows of $\lambda$-displacement convex functionals in $\mathcal P_2(\rr^n)$, see~\cite{ambrosio-gigli-savare}). 

Here also $H(\mu_{t} \vert \nu)$ and $I(\mu_{t} \vert \nu)$ get instantaneously finite  for $\mu_0 \in \mathcal P_2(\rr^n)$ if $\nabla^S A$ is uniformly bounded from below, as can be seen by adapting the proofs of the Theorem in~\cite{ov01} and Lemma~\ref{cor-2} below.

Let us also notice that the coupled conditions $\int |F|^4d\nu<\infty, \int \mu_0^2 e^{V} <  \infty$ can be modified by using the H\"older or the Young inequality instead of the Cauchy-Schwarz inequality, and for instance be replaced by $\int e^{F^2} d\nu < \infty, H(\mu_0 \vert \nu) < \infty$.
\end{erem}

\begin{ecor}
\label{cor-cvWJ}
Let $d\nu=e^{-V}dx \in \mathcal P_2(\rr^n)$ with $\nabla \cdot(e^{-V} F) =0$ and make the same hypotheses as in Theorem~\ref{thm-derivation}.
Assume moreover the existence of a constant $C>0$ such that
\begin{equation}
\label{eq-def-wj}
W_2^2(\mu,\nu)\leq \frac1C\int (x-\n\psi)\cdot \PAR{\nabla\log\mu+A} \,d\mu 
\end{equation}
for all probability densities $\mu$, where ${\nabla\psi}\#\mu=\nu$. Then 
\begin{equation}\label{eq:contr2}
W_2(\mu_t,\nu)\leq e^{-Ct}W_2(\mu_0,\nu),  \qquad t \geq 0
\end{equation}
for any solution $(\mu_t)_t$ to~\eqref{eq-fp} starting from a probability density $\mu_0$ as in Theorem~\ref{thm-derivation}.
\end{ecor}

\begin{eproof}
It is a consequence of Theorem~\ref{thm-derivation} since the map $t\mapsto W_2(\mu_t,\nu)$ is absolutely continuous.
\end{eproof}

\medskip

We saw in the introduction that the contraction property~\eqref{eq-contraction-wasserstein} between all solutions, whence the uniform exponential convergence estimate~\eqref{eq-wasserstein-decay}-\eqref{eq:contr2}, holds if $A = \nabla V$ with $\nabla^2 V (x) \geq C {\rm Id}_n$ for all $x$, or more generally if $\nabla^S A(x) \geq  C {\rm Id}_n$. 

Let us now give a first simple and weaker criterion ensuring the condition~\eqref{eq-def-wj} in Corollary~\ref{cor-cvWJ}, whence the uniform exponential convergence~\eqref{eq:contr2} of the solutions to $\nu$. 

For that purpose, recall that a measure $\nu$ is said to satisfy a (transportation) Talagrand inequality with constant $C>0$, denoted $WH(C)$, if
\begin{equation}\label{WH}
W_2(\nu,\mu)\le \sqrt{\frac2C \, H(\mu \vert \nu)}
\end{equation}
for all measure $\mu$ absolutely continuous with respect to $\nu$ (see~\cite{ov00} for instance). Then :
 
\begin{eprop}\label{prop:AGSWH}
Assume that the measure $d\nu = e^{-V} dx \in \mathcal P_2(\rr^n)$ satisfies a WH(c) inequality and that $\nabla^2 V(x) \geq \lambda_1 {\rm Id}_n$, $\nabla^S F(x) \geq \lambda_2 {\rm Id}_n$ with $\lambda_1, \lambda_2 \in \rr$ and all $x$. Then it satisfies~\eqref{eq-def-wj} with the constant $C =  (c+ \lambda_1 + 2 \lambda_2)/2)$ inequality if $c > - \lambda_1 - 2 \lambda_2$. 
\end{eprop}

In particular, when $A = \nabla V$ and $\lambda_1>0$, then $\nu = e^{-V}$ satisfies a WH($\lambda_1$) inequality (see~\cite{ov00}), so that~\eqref{eq-def-wj} holds with the constant $\lambda_1$, as observed above (see also Lemma~\ref{prop:WJunifcvx} below for the non-gradient case). But above all it allows for larger classes of measures $\nu$ satisfying a $WH$ inequality, as described in~\cite{gozlanleonard}, including for example potentials $V$ which are the sum of a uniformly convex and of a bounded function. 

\smallskip

 \begin{eproof}
By Lemma~\ref{cor-2} below and assumptions,
 $$
- \int (x-\n\psi)\cdot(\n\log{\mu}+A)d\mu
+ \Big(\frac{\lambda_1}{2} + \lambda_2 \Big) \, W_2^2(\nu,\mu) 
\leq
- H(\mu \vert \nu) \leq - \frac{C}{2} W_2^2(\nu,\mu) 
$$
for any probability density $\mu$ with $\nabla \psi \# \mu = \nu$. This concludes the argument.
\end{eproof}

\begin{elem}
\label{cor-2}
Let $d\nu=e^{-V}dx \in \mathcal P_2(\rr^n)$ and $F$ be a vector field such that $\n~\cdot~(e^{-V}F)=0$. If $\nabla^2 V  \geq \lambda_1 {\rm Id}_n$ and $\nabla^S F  \geq \lambda_2 {\rm Id}_n$
for some $\lambda_1,\lambda_2 \in \rr$,  uniformly in $\dR^n$, then 
\begin{equation}
\label{eq-anous}
H(\mu|\nu)+\big(\frac{\lambda_1}{2}+\lambda_2\big)W_2^2(\nu,\mu)\leq \int (x-\n\psi)\cdot(\n\log{\mu}+A)d\mu
\end{equation}
for every probability density $\mu$, and with ${\nabla\psi}\#\mu=\nu$. 
\end{elem}

\begin{eproof}
We follow the proof of Theorem~1 in~\cite{cordero}.  Let $\mu$ be a probability on $\dR^n$ with a smooth positive density $f$ with respect to $\nu$. If $\n\psi\#\mu=\nu$ then, by change of variables,
$$
fe^{-V} = 
e^{-V(\n\psi)}\det (\n^2\psi). 
$$
Then 
\begin{eqnarray*}
\int f\log f d\nu
&=&
\int [V-V(\n\psi)+\log\det (\n^2\psi)]fd\nu
\leq
 \int [V-V(\n\psi)+\Delta(\psi-\frac{|x|^2}{2})]fd\nu\\
&\leq&
\int [V-V(\n\psi)+\n V\cdot(\n\psi-x)]fd\nu-\int(\n\psi-x)\cdot\n fd\nu 
\end{eqnarray*}
by convexity and integration by parts. Here $\Delta \psi$ is the Alexandrov Laplacian of the convex function $\psi$, which is smaller than its distributional Laplacian. Moreover  
$$
\int (x-\n\psi)\cdot(\n\log{\mu}+A)d\mu =\int (x-\n\psi)\cdot\n f\,d\nu+\int (x-\n\psi)\cdot F\,f\,d\nu
$$
and
$$
\int (x-\n\psi)\cdot F(\n\psi)\,d\mu= \int (\nabla \psi^* -x) \cdot  F \, d\nu = - \int (\psi^* - \frac{\vert x \vert^2}{2}) \nabla \cdot(e^{-V} F) =0
$$
since $\n \psi \# \mu=\nu$ and $\n \cdot (e^{-V}F)=0$. Hence
\begin{eqnarray*}
H(\mu|\nu) 
= \int f\log f d\nu 
&\leq&
 \int [V-V(\n\psi)+\n V\cdot(\n\psi-x)-(F- F(\n \psi)) \cdot(x-\n\psi)]d\mu
 \\
 && + \int (x-\n\psi)\cdot(\n\log{\mu}+A)d\mu
. 
\end{eqnarray*}

Now, by a Taylor expansion,
$$
V-V(\n\psi)+\n V\cdot(\n\psi-x)
= - \int_0^1(1-t)(\n\psi -x) \cdot [\n^2V(x+t(\n\psi-x))(\n\psi - x )]dt\leq-\frac{\lambda_1}{2} |\n\psi -x|^2
$$
and 
$$
- (F-F(\n\psi))\cdot(x-\n\psi)=-\int_0^1(\n\psi -x)\cdot [\n^SF(x+t(\n\psi-x))(\n\psi -x)]dt \leq-{\lambda_2}|\n\psi -x|^2. 
$$
This concludes the argument by combining the two expressions. 
\end{eproof}

\begin{erem}\label{rem:AGSOV}
When $F=0$, then inequality~\eqref{eq-anous} has been derived in~\cite{ov00} in the proof of the HWI inequality
$$
H(\mu|\nu) \leq W_2(\nu,\mu) \, \sqrt{I(\mu \vert \nu)} - \frac{\lambda_1}{2}W_2^2(\nu,\mu),
$$ 
where $I(\mu \vert \nu)$ is the Fisher information of $\mu$ with respect to $\nu$, defined in~\eqref{eq-fisher}. It implies the HWI inequality since
 $$
 \int (x-\n\psi)\cdot\n\log \frac{\mu}{\nu}\,d\mu \leq \sqrt{\int |x-\n\psi|^2\,d\mu} \; \sqrt{\int \ABS{\n\log \frac{\mu}{\nu}}^2\,d\mu}
= W_2(\nu,\mu) \, \sqrt{I(\mu \vert \nu)}
$$
by the Cauchy-Schwarz inequality; here again ${\nabla\psi} \#\mu=\nu.$

Moreover, again for $F=0$, inequality~\eqref{eq-anous} appears in~\cite{ambrosio-gigli-savare} as a fundamental inequality in the general theory of gradient flows, see~\cite[Th.~4.0.4]{ambrosio-gigli-savare} for instance.
\end{erem}

\begin{erem}
In this work we focus on the estimate~\eqref{eq-wasserstein-decay} in the Euclidean Wasserstein distance and give simple necessary and sufficient conditions (weaker than strictly positive curvature) on the drift for~\eqref{eq-wasserstein-decay} to hold for any initial condition $\mu_0$. 

Let us stress  that in our study, it is important that there is no (larger than 1) multiplicative constant on the right-hand side of (\ref{eq-wasserstein-decay}). Indeed, there are various ways to get convergence result of the form
\begin{equation}\label{K}
W_2(\mu_t,\nu)\le K\,e^{-Ct}\,W_2(\mu_0,\nu)
\end{equation}
for a constant $K$ larger than 1. Let us mention two different approaches.
\begin{enumerate}
\item Suppose that $\nu$ satisfies a logarithmic Sobolev inequality with constant $C$, that is
\begin{equation}\label{LSI}
H(f \nu \vert \nu) \leq \frac{1}{2C} I(f\nu \vert \nu)
\end{equation}
for all probability densities $f$ with respect to $\nu$. This inequality can be proved in infinite negative curvature cases
and is equivalent to the exponential decay of the entropy
$$
H(\mu_t \vert \nu) \le e^{-2C(t-t_0)} H(\mu_{t_0} \vert \nu).
$$
Recall then that such a logarithmic Sobolev inequality implies the Talagrand inequality~\eqref{WH} with the same constant $C$ (see for example \cite{ov00}).
 Hence 
$$
W_2(\mu_t,\nu)\le K(C,t_0)\,e^{-Ct}\,\sqrt{H(\mu_{t_0} \vert \nu)} \le \tilde K(V,C,t_0) \,e^{-C t}\, W_2(\nu,\mu_0)
$$
for all $t$. The last inequality follows from a regularization argument derived from a Harnack type inequality under regularity assumptions on $V$ (see \cite{wang-book}).
\item Another approach relies on the study of the contraction in a Wasserstein distance for a twisted metric, equivalent to the Euclidean one, so that such a contraction result will lead to convergence in the Euclidean Wasserstein distance as in~\eqref{K}, with a $K>1$. This has been successfully done for the kinetic Fokker-Planck equation in a perturbation of the Gaussian case (infinite curvature case) in \cite{bgm} using the simplest coupling (same Brownian motion for the two different dynamics, as in the introduction). Recently, A.~Eberle \cite{ebe2011} has used reflection coupling to establish contraction results in a twisted metric for a reversible Fokker-Planck equation under lower negative curvature and sufficient quadratic growth condition at infinity.
\end{enumerate}
\end{erem}

 \section{The $WJ$ inequality }
 \label{sec-wj}

In this section we derive a functional inequality ensuring the uniform exponential convergence~\eqref{eq:contr2} of the solutions to the steady state $\nu = e^{-V}$, give practical criteria for it and its connections with classical functional inequalities as the Poincar\'e or logarithmic Sobolev inequalities. 

\subsection{Definition of the inequality}

As in Theorem~\ref{thm-derivation}, let us assume that $V,F$ are such that $\int |F|^4d\nu<\infty$ with $d\nu = e^{-V} dx$. 
If $\mu_t$ is a solution of~\eqref{eq-fp} with initial condition having a smooth density $\mu_0$ such that  
$\int \mu_0^2 e^{V}dx<\infty$, we saw that the map $t\mapsto W_2(\mu_t,\nu)$ is absolutely continuous and for  almost every $t \geq 0$
$$
\frac12 \frac{d}{d t}W_2^2(\mu_t,\nu)=\int \PAR{\nabla \psi_t (y) -y}\cdot \PAR{\nabla \log\mu_t (y)+A(y)}d\mu_t (y)
$$
where for every $t\geq 0$, ${\nabla \psi_t}\#{\mu_t}=\nu$.  In fact, since $d\nu = e^{-V} dx$ is a stationary solution of~\eqref{eq-fp}, and since $\int \vert F \vert^2 d\nu < \infty$, then, again by~\cite[Th.~23.9]{villani-book1}, 
\begin{eqnarray*}
&&\frac12 \frac{d}{d t}W_2^2(\mu_t,\nu)
\\
&&=
\int \PAR{\nabla \psi_t (y)-y}\cdot \PAR{\nabla \log\mu_t (y)+A(y)} d\mu_t (y)+ \int \PAR{\nabla \phi_t(x) - x}\cdot \PAR{\nabla \log\nu (x)+A(x)} d\nu (x).
\\
&&=
\int \PAR{\nabla \psi_t (y)-y}\cdot \PAR{\nabla \mu_t (y)+A(y) \mu_t (y)} dy + \int \PAR{\nabla \phi_t(x) - x}\cdot \PAR{\nabla \nu (x)+A(x) \nu (x)} dx.
\end{eqnarray*}
Here ${\nabla \phi_t}\#{\nu}=\mu_t$, so that $\psi_t = \phi_t^*$.

Then one can perform a ``weak" integration by parts as in~\cite[Th. 1.5]{lisini} and use the push-forward property to bound from above the right-hand side by
$$
- \int_{\rr^n} \Big( \Delta\phi_t (x)+\Delta\phi_t^*(\nabla\phi_t(x))-2n +(A(\nabla\phi_t(x))- A(x)) \cdot (\nabla\phi_t(x)-x) \Big) \,d\nu(x).
$$
Here $\Delta\phi$ is the Alexandrov Laplacian of a convex map $\varphi$ on $\rr^n$. 

Observe now that for $t>0$ both $\nu$ and $\mu_t$ belong to the set $\mathcal P_{2,c} (\rr^n)$ of measures of $\mathcal P_2(\rr^n)$ with $\mathcal C^1$ positive densities on $\rr^n$. In particular, a measure $\rho$  in $\mathcal P_{2,c} (\rr^n)$ has a density which is $\mathcal C^{0, \alpha}$ and bounded from above and from below by a positive constant on any ball of $\rr^n$.  Then Caffarelli's regularity results (see~\cite{caffarelli92}) also apply in the case of two measures $\rho_1$ and $\rho_2$ in $\mathcal P_{2,c} (\rr^n)$, and ensure that both convex functions $\phi$ and $\phi^*$, where $\nabla \phi \# \rho_1 = \rho_2$ and 
$\nabla \phi^* \# \rho_2 = \rho_1$, are $\mathcal C^2$ and strictly convex. 

In particular here the convex functions $\phi_t$ and $\phi_t^*$ are $C^2$ and strictly convex, and $\Delta\phi_t$ and $\Delta\phi_t^*$ are the usual Laplacians; moreover for almost every $t \geq 0$
\begin{multline*}
\frac12 \frac{d}{d t}W_2^2(\mu_t,\nu)
\\
\leq
- \int_{\rr^n} \Big( \Delta\phi_t (x)+\Delta\phi_t^*(\nabla\phi_t(x))-2n +(A(\nabla\phi_t(x))- A(x)) \cdot (\nabla\phi_t(x)-x) \Big) \,d\nu(x).
\end{multline*}

This motivates the following definition:

\begin{edefi}
We say that the couple $(\nu,A)$, where $\nu$ belongs to $\mathcal P_{2,c} (\rr^n)$ and $A$ is a $\mathcal C^1$ vector field,  satisfies a $WJ$ inequality with constant $C>0$ if 
\begin{equation}
\label{eq-defn}
 W_2(\nu,\mu)\leq\sqrt{\frac1C \,J(\mu|(\nu,A))} 
\end{equation}
for every $\mu \in \mathcal P_{2,c} (\rr^n)$; here
$$
J(\mu|(\nu, A))= \int  \big[ \Delta\phi+ \Delta\phi^*(\nabla\phi) - 2n+ (A(\nabla\phi) - A) \cdot (\nabla\phi-x) \big] d\nu
$$
where ${\nabla\phi}\#\nu=\mu$. We implicitly assume in the definition that $J(\mu|(\nu, A))$ is well defined and non-negative.

For simplicity, if $d\nu=e^{-V}dx$ and $A=\n V$, or equivalently $F=0$,  then $J(\mu|(\nu, A))$ is denoted  $J(\mu|\nu)$ and we say that the probability measure  $\nu$ satisfies a $WJ$ inequality.  
\end{edefi}

\bigskip

This definition is general, and does not assume that $\nu$ is invariant with respect to the Fokker-Planck equation driven by $A$; when it is the case, that is, when $d\nu = e^{-V} dx$ and $\nabla \cdot (e^{-V} F)) =0$, then as in Corollary~\ref{cor-cvWJ}, the $WJ$ inequality governs the uniform exponential convergence of solutions to~\eqref{eq-fp1} towards the equilibrium $\nu$, according to~\eqref{eq:contr2}.

\bigskip
\subsection{Sufficient conditions}

We begin with the following simple but key observation :

\begin{elem}
\label{lem-commode}
If $\phi$ is a $\mathcal C^2$ strictly convex function on $\rr^n$ then
$$
 \Delta \phi (x) + \Delta\phi^*(\nabla\phi (x)) - 2n \geq0
$$ 
for all $x$ such that the Hessian matrix $\nabla^2 \phi (x)$ at $x$ is positive, and is $0$ if and only if $\nabla^2 \varphi(x)$ is the identity matrix.
\end{elem}

\begin{eproof}
 Given $x \in \rr^n$ we write $\nabla^2\phi (x)$  as $O\,D\,O^*$ where $O$ is orthonormal, $D=diag(d_1,..,d_n)$ and $d_i$ are the positive eigenvalues of $\nabla^2\phi (x)$. 

Observe that  $\nabla \phi^*(\nabla \phi(x))=x$, and then 
$$
\nabla^2\phi^*(\nabla \phi(x))\nabla^2\phi(x)={\rm Id}_n.
$$
 This leads to
 $$\nabla^2\phi^*(\nabla\phi (x))=(\nabla^2\phi(x))^{-1}=O\,D^{-1}\,O^*.$$
Then
 $$
 \Delta\phi (x) +  \Delta\phi^*(\nabla\phi (x)) - 2n= \sum_{i=1}^n d_i + \sum_{i=1}^n \frac{1}{d_i} - 2n = \sum_{i=1}^n\left(\sqrt{d_i}-\frac1{\sqrt{d_i}}\right)^2 \geq0,
 $$
 with equality if and only if the $d_i$ are all equal to $1$.
\end{eproof}

If $A$ is {\it monotone}, that is, if
$$
(A(x) - A(y)) \cdot (x-y) \geq 0
$$
for all $x,y$, then by Lemma~\ref{lem-commode} the quantity $J(\mu \vert (\nu, A))$ is non-negative for all $\mu = {\nabla\phi}\# \nu$. This is not always the case, as pointed out to us by B. Han (see~\cite{han}). Observe similarly that along the evolution of the Fokker-Planck equation, the dissipation of the relative entropy to the steady state, and more generally of relative $\phi$-entropies with $\phi$ convex, is non-negative; this is however not always the case for the Fisher information, as observed by B. Helffer (see~\cite{logsob}).

\bigskip


Lemma~\ref{lem-commode} has the following straightforward consequence :

\begin{elem}\label{prop:WJunifcvx}
If $\nu$ is in $\mathcal P_{2,c}(\rr^n)$ and $A$ is such that 
\begin{equation}\label{eq:Aunif}
\nabla^S A \geq C \, {\rm Id}_n
\end{equation}
with $C>0$, uniformly on $\rr^n$, then $(\nu, A)$ satisfies a $WJ(C)$ inequality.
\end{elem}  

\smallskip

This is natural since the contraction property~\eqref{eq-contraction-wasserstein} between any solutions, and not only the convergence estimate~\eqref{eq-wasserstein-decay}, holds in this uniformly monotone situation, as observed in the introduction.

\smallskip
In particular the standard Gaussian measure $\gamma$ on $\rr^n$ satisfies a $WJ$ inequality with constant $1$ and  the constant $1$ is optimal. Observe indeed that
$$
J(\mu| \gamma) = \int  \PAR{\Delta\phi (x) +\Delta\phi^*(\nabla\phi(x)) - 2n} d\gamma(x) + W_2^2(\mu, \gamma)
$$
for all $\mu = \nabla \varphi \# \gamma$. Hence it is always larger than $W_2^2(\mu, \gamma)$ by Lemma~\ref{lem-commode}; moreover it is equal to $W_2^2(\mu, \gamma)$ if and only if the non-negative term $\Delta\phi (x) +\Delta\phi^*(\nabla\phi(x)) - 2n$ is $0$ for almost every $x$, that is, if and only if $\nabla^2 \varphi(x) = {\rm Id}_n$, by Lemma~\ref{lem-commode}, that is, if and only if $\mu$ is a translation of~$\gamma$.

\bigskip

For uniformly convex potentials $V$, or more generally under~\eqref{eq:Aunif}, the $WJ$ inequality for $(\nu, A)$ is obtained without using the non-negative contribution  $\Delta\phi (x) +\Delta\phi^*(\nabla\phi(x)) - 2n$ in $J$, which stems from the diffusion term. 

Proposition~\ref{prop:AGSWH} gave a first way of taking advantage of  the diffusion term to consider non-uniformly convex cases and even non-convex cases. However, for a non-gradient drift, there is a strong assumption on the measure $\nu$ in Proposition~\ref{prop:AGSWH}, which is not always easy to be checked since $\nu$ may not be explicit. We can replace it by another criterion, which asks for weaker assumptions on $\nu$, for instance:

\begin{eprop}\label{Amonotone}
Let $A$ be a $C^1$ monotone map from $\rr^n$ to $\rr^n$ for which there exist two constants $R \geq 0$ and $K>0$  such that
$$
\nabla^S A(x)\geq K
$$
 for all $\vert x \vert \geq R$, and let $d\nu = e^{-V}$ be a probability measure on $\rr^n$, with $V$ a $\mathcal C^1$ potential.
 
Then $(\nu, A)$ satisfies a $WJ$ inequality with constant  $C=C(V,R,K)$.
\end{eprop}

\begin{erem}
The constant $C$ given by the proof depends on $V$ only through its minima and maxima on the ball of center $0$ and radius $3R$. Observe that the proof requires only $V$ to be bounded on this ball, and that any ball of center $0$ and radius $>R$ would work. 
\end{erem}

The proof consists in overcoming the lack of convexity near the origin by using the diffusion term.  It will be given at the end of the section.

\medskip

Let us see  the influence of the diffusion term on the rate of convergence to equilibrium on a simple example, for instance for the potential $V(x) = x^4$ on $\rr$ and $F=0$. By Proposition~\ref{prop:AGSWH} or~\ref{Amonotone}, the measure $e^{-V}$ satisfies the condition~\eqref{eq-def-wj} with a constant $C>0$, whence solutions $\mu_t$ to the Fokker-Planck equation~\eqref{eq-fp1} converge exponentially fast to it, according to $W_2(\mu_t, e^{-V}) \leq  e^{-Ct}W_2(\mu_t, e^{-V})$.  On the other hand, without diffusion, the solution at time $t$ to $\partial_t \mu_t = \nabla \cdot(\mu_t \nabla V)$ is the distribution of the points $x(t)$ initially at $x(0)$ drawn according to $\mu_0$ and evolving according to $x'(t) = - x^3$. This solves into $x(t)^2 = x(0)^2/(1+2t x(0)^2)$, so that the solution $\mu_t$ converges to the unique steady state $\delta_0$ according to $$
W_2^2(\mu_t, \delta_0) = \int \frac{x^2}{1+2t x^2} \, d\mu_0(x) \sim \frac1t
$$
for large $t$.

\begin{erem}\label{rem:sturmvonrenesse}

If $(\mu_t)_t$ and $(\nu_t)_t$ are two solutions to~\eqref{eq-fp}, a formal adaptation of the above computation gives
$$
\frac{1}{2} \frac{d}{dt} W_2^2(\mu_t, \nu_t) 
= - \int  \big[ \Delta\phi_t+ \Delta\phi_t^*(\nabla\phi_t) - 2n+ (A(\nabla\phi_t) - A) \cdot (\nabla\phi_t-x) \big] d\mu_t
$$ 
if $\nu_t = \n \varphi_t \# \mu_t$. With this in hand one can recover the equivalence between the following three assertions, due to K.-T. Sturm and M. von Renesse (see~\cite{sturm-vonrenesse}  and \cite[Th. 5.6.1]{wang-book})~:

1) For all initial conditions $\mu_0$ and $\nu_0$ in $\mathcal P_2(\rr^n)$, for all $t\geq0$,
$$
W_2(\mu_t,\nu_t)\leq e^{- Ct}W_2(\mu_0,\nu_0),
$$
where $\mu_t$ (resp. $\nu_t$) are solutions of~\eqref{eq-fp1} starting from $\mu_0$ (resp. $\nu_0$). 

1') For all $x,y\in\dR^n$ and all $t\geq0$,
$$
W_2(\mu_t,\nu_t)\leq e^{- Ct}|x-y|,
$$
where $\mu_t$ (resp. $\nu_t$) are solutions of~\eqref{eq-fp1} starting from $\delta_x$ (resp. $\delta_y$).

2) For all $x,y\in\dR^n$ the vector field  $A$ satisfies 
$$
(A(y) - A(x)) \cdot (y-x) \geq  C \vert y-x \vert^2.
$$

Indeed time-differentiating $1')$ at $t=0$ implies $2)$, and $2)$ implies $1)$ by time-integration and Lemma~\ref{lem-commode}.
\end{erem}

 \subsection{Tensorization and perturbation}

Fundamental properties of functional inequalities lie in the range of stability: non dependence on the dimension, which enables to consider problems in infinite dimension, and stability by perturbation, which enables to reach more general potentials.

The following two results are important to extend the practical conditions we just derived. The first one concerns the tenzorization : namely, the product of measures satisfying a $WJ$ inequality also satisfies a $WJ$ inequality.

\begin{eprop}[Tensorization]
\label{tenso}
Suppose that the measures and drifts $(\nu_i, A_i)_{1\le i\le n}$ satisfy a $WJ(C_i)$ inequality on $\R^{n_i}$ respectively. Then $(\displaystyle \otimes_{i=1}^n\nu_i, A)$ with $A(x) = (A_i(x_i))_{1 \leq i \leq n}$ for $x = (x_i)_i$ on the product space satisfies a $WJ$ inequality with constant $\min_i C_i$.
\end{eprop}

\begin{eproof}
Let us assume for simplicity of notation that $n_i=1$ for all $i$, and let us denote $\displaystyle d\nu^n (x) =\otimes_{i=1}^n d\nu_i(x_i) \, dx_i $. For $x \in \rr^n$ we let  $\hat x_i \in \rr^{n-1}$ have the same coordinate than $x$, but the i-th coordinate $x_i$, which is removed. 

Let now $\phi$ be a $\mathcal C^2$ strictly convex function on $\rr^n$. Noticing that all its restrictions $x_i \mapsto \varphi(\hat x_i, x_i)$ are also $\mathcal C^2$ strictly convex functions on $\rr$, and using the $WJ$ inequality for each $\nu_i$ we get 
\begin{eqnarray*}
\int_{\rr^n} \! |\nabla\phi(x)-x|^2d\nu^n (x)
&=&
\sum_{i=1}^n\int_{\rr^{n-1}} \displaystyle \otimes_{j\not=i} d\nu_j(\hat x_i) \int_{\rr} |\partial_i\phi(x)-x_i|^2d\nu_i(x_i)\\
&\le&\! \! 
\frac1{\min_i C_i}\sum_{i=1}^n\int \! \otimes_{j\not=i}d\nu_i(\hat x_i)\int \! (A_i(\partial_i\phi(x))-A_i(x_i))(\partial_i\phi(x)-x_i)d\nu_i (x_i)\\
&&
\qquad+\frac1{\min_i C_i}\sum_{i=1}^n\int \! \otimes_{j\not=i} d\nu_j(\hat x_i )\int \! \left(\partial^2_{ii}\phi(x)+\frac{1}{\partial^2_{ii}\phi(x)}-2\right)d\nu_i (x_i) \\
&\le&
\frac1{\min_i C_i}\int_{\rr^n} \sum_{i=1}^n(A_i(\partial_i\phi(x))- A_i(x_i))(\partial_i\phi(x)-x_i) d\nu^n(x) \\
&&\qquad\qquad+\frac1{\min_i C_i}\int_{\rr^n} \sum_{i=1}^n \left( \partial^2_{ii}\phi(x)+\frac{1}{\partial^2_{ii}\phi(x)}-2 \right)d\nu^n(x).
\end{eqnarray*}

Now, in the first term,
$$
\sum_{i=1}^n(A_i(\partial_i\phi(x))-A_i(x_i))(\partial_i\phi(x)-x_i)=(A(\nabla\phi(x))-A(x)).(\nabla\phi(x)-x).
$$

In the second term we fix $x \in \rr^n$ and, in the notation of Lemma~\ref{lem-commode}, we write $\nabla^2\phi(x)= O D O^*$ where $O$ is orthonormal and $D = diag(d_1, \dots, d_n)$. Then
$$
\sum_{i=1}^n\partial^2_{ii}\phi (x) =tr(\nabla^2\phi (x))=\sum_{i=1}^nd_i.
$$
Moreover $\partial_{ii}\phi(x)=\sum_{j=1}^nO_{ij}^2d_j$ with $\sum_{i=1}^nO_{ij}^2=1$, and $x\mapsto x^{-1}$ is convex on $\{x>0\}$, so by the Jensen inequality
$$
\sum_{i=1}^n\frac{1}{\partial^2_{ii}\phi (x)}
=
\sum_{i=1}^n\frac1{\sum_{j=1}^nO_{ij}^2d_j}
\le
\sum_{i=1}^n \; \sum_{j=1}^nO_{ij}^2\frac1{d_j} 
=
\sum_{j=1}^n \; \sum_{i=1}^nO_{ij}^2\frac1{d_j}
=
\sum_{j=1}^n\frac1{d_j}
$$
since also $\sum_{j=1}^nO_{ji}^2=1$. Hence
$$
\sum_{i=1}^n \left( \partial^2_{ii}\phi(x)+\frac{1}{\partial^2_{ii}\phi(x)}-2 \right) \leq \sum_{i=1}^n \left( d_i +\frac{1}{d_i}-2 \right) = \Delta \phi (x) + \Delta\phi^*(\nabla\phi (x)) -2n
$$
as in the proof of Lemma~\ref{lem-commode}. This concludes the proof.
\end{eproof}

Let us come back to the PDE motivation of the WJ inequality : letting $d\nu_i(x_i) = e^{-V_i(x_i)} dx_i$ for each $i$, then $\nabla \cdot ((A-\nabla V)(x) e^{-V(x)}) =0$ on the product space with $V(x) = \sum_{i=1}^n V_i(x_i)$ for $x = (x_i)_i$ as soon as $\nabla \cdot ((A_i - \nabla V_i)(x_i) e^{-V_i(x_i)}) =0$ on $\rr^{n_i}$ for each~$i$. Hence $e^{-V} dx = \displaystyle \otimes_{i=1}^n e^{-V_i(x_i)} dx_i$ is indeed a stationary measure of the corresponding PDE on the product space if so is each $e^{-V_i(x_i)} dx_i$ on $\rr^{n_i}$.

\bigskip

The second result is about the perturbation of the measure $\mu$. For classical functional inequalities, such as Poincar\'e inequality, $WH$ or  logarithmic Sobolev inequality, perturbations by bounded potentials are allowed (see for example~\cite{logsob,gozlanleonard}). Here we have to be more restrictive, not only on the perturbation term but also on the initial measure satisfying a $WJ$ inequality.

\begin{eprop}[Perturbation]
\label{perturb}
Suppose that the measure and drift $(\nu, A)$ satisfy a $WJ(C)$ inequality and that for an $\alpha \leq 0$
\begin{enumerate}
\item $(A(y) - A(x)) \cdot (y-x) \geq \alpha \vert y-x \vert^2$ for  all $x,y$.
\suspend{enumerate}
Consider a map $T$ on $\rr^n$ such that  $e^{-T}d\nu$ is a probability measure and  for a $K \geq 0$
\resume{enumerate}
\item $|T(x)| \le K$ for all $x$,
\suspend{enumerate}
and a map $B$ from $\rr^n$ to $\rr^n$ such that for a $\beta \in \rr$
\resume{enumerate}
\item $(B(y) - B(x)) \cdot (y-x) \geq \beta \vert x-y \vert^2$ for all $x,y$.
\end{enumerate}
If $ - \beta e^{2K} - \alpha (e^{2K} -1)< C$, then $(e^{-T} \nu, A+B)$ satisfies a $WJ$ inequality with constant $Ce^{-2K} +\beta + \alpha (1 - e^{-2K})$.
\end{eprop}

\begin{eproof}
Let $\tilde{\nu} = e^{-T} \nu$, and let $\phi$ be a $C^2$ strictly convex map on $\rr^n$. Then 
\begin{eqnarray}
\int|\nabla\phi(x)-x|^2d\tilde\nu (x)
&\stackrel{ii.}{\le}&
e^K\int|\nabla\phi(x)-x|^2d\nu \nonumber\\
&\stackrel{WJ}{\le}&
\frac{e^K}C\int(A(\nabla\phi)-A) \cdot (\nabla\phi-x)d\nu\nonumber\\
&&+\frac{e^K}C\int(\Delta\phi + \Delta\phi^*(\nabla\phi)-2n)d\nu. \label{2terms}
\end{eqnarray}
Since $\Delta\phi(x) + \Delta\phi^*(\nabla\phi(x))-2n\ge0$ by Lemma~\ref{lem-commode}, the second integral on the right-hand side of~\eqref{2terms} is bounded by
$$
 \frac{e^{2K}}{C} \int(\Delta\phi + \Delta\phi^*(\nabla\phi) -2n)d\tilde\nu
$$
by $ii.$ Moreover, by $i.$ and $ii.$, we write the first integral on the right-hand side of~\eqref{2terms} as
\begin{multline*}
 \int \big[ (A(\nabla\phi)-A) \cdot (\nabla\phi-x) - \alpha \vert \nabla \phi - x \vert^2 \big]d\nu + \alpha  \int \vert \nabla \phi - x \vert^2d\nu\\
 \leq e^K \int(A(\nabla\phi)-A) \cdot (\nabla\phi-x) \, d\tilde\nu - \alpha  (e^K - e^{-K}) \int \vert \nabla \phi - x \vert^2 \, d\tilde\nu.\\
\end{multline*}
Then, by $iii.$, we bound the first integral on the above right-hand side by
$$
\int((A+B)(\nabla\phi)-(A+B)(x)).(\nabla\phi-x)d\tilde\nu
- \beta \int \vert \nabla \phi - x \vert^2 \, d\tilde\nu.
$$

This concludes the proof by collecting all terms and using the positivity conditions on the coefficients.
\end{eproof}

Typically $A=\nabla V$ with $\nu=e^{-V}$ and the bounded perturbation is given by $B=\nabla T$.  
Note also that one can adapt the proof above to give a variant of this result for $\alpha > 0$.

 \subsection{Necessary conditions}

We now compare the $WJ$ inequality for a measure $\nu = e^{-V}$ and a drift $A $ with more classical inequalities.
 
 \smallskip
 
We first prove that a $WJ$ inequality implies a Poincar\'e inequality:
 
 \begin{eprop}\label{WJ-PI}
 If $(\nu, A)$ satisfies a $WJ(C)$ inequality then $\nu$ satisfies a Poincar\'e inequality with the same constant $C$, that is, for every smooth function $f$ 
 $$
\int \PAR{f-\int fd\mu}^2d\nu\leq\frac1C \int |\nabla f|^2d\nu. 
 $$
 \end{eprop}
 \begin{eproof}
 Let $f$ be a smooth map on $\rr^n$ and $\varphi$ be defined by
$$
 \varphi(x)=\frac{|x|^2}2+\varepsilon f(x)
 $$
 for small $\varepsilon$. Then for all $x$ the Hessian matrices of $\varphi$ and $f$ and their respective eigenvalues $d_i$ and $f_i$ for $1 \leq i \leq n$ satisfy
$$
\nabla^2\phi(x)={\rm{Id}_n}+\varepsilon \,\nabla^2 f(x),\qquad \qquad d_i=1+\varepsilon\, f_i.
$$
Hence, as in Lemma~\ref{lem-commode}, 
 \begin{eqnarray*}
\Delta\phi^*(\nabla\phi(x)) + \Delta \phi (x) -2n
=
 \sum_{i=1}^n \left(\frac 1{d_i}+d_i-2\right)
&=&\sum_{i=1}^n \left(\frac1{1+\varepsilon f_i}+1+\varepsilon f_i -2\right)\\
&=&\varepsilon^2\sum_{i=1}^n  f_i^2+o(\varepsilon ^2)
=\varepsilon^2 \| \nabla^2 f\|_{HS}+o(\varepsilon^2).
 \end{eqnarray*}
Moreover
 $$\nabla V(\nabla\phi(x))-\nabla V(x)=\varepsilon \, \nabla^2V(x) \nabla f(x)+o(\varepsilon).$$
Hence, for this map $\varphi$, the $WJ$ inequality now reads
 $$
 \varepsilon^2\int\left[  \|\nabla^2 f(x)\|_{HS}+\nabla f(x)\cdot\nabla^2V(x) \nabla f(x)\right]d\nu(x) +o(\varepsilon^2)\ge \varepsilon^2 \, C\,\int|\nabla f(x)|^2d\nu
 $$
where $\|M\|_{HS}$ is the Hilbert-Schmidt norm of a matrix $M$. Letting $\varepsilon\to0$, we recover the well-known integral $\Gamma_2$ criterion (see for example \cite[Prop. 5.5.4]{logsob}), which is equivalent to the Poincar\'e inequality with constant $C$.
 \end{eproof}

We now turn to the $WI$ inequality in the particular case $A=\nabla V$.
 
An inequality looking like $WJ$ has been introduced in~\cite{ov00} and studied in~\cite{glwy,glww} for its equivalence to deviation inequalities for integral functional of Markov processes: thus it has high practical interest. We say that a probability  measure $\nu$ satisfies a $WI$ inequality with constant $C>0$ (called $LSI+T(C)$ in~\cite{ov00}) if for every probability measure $\mu$ absolutely continuous with respect to $\nu$  
$$
W_2(\nu,\mu)\leq \frac1C \sqrt{I(\mu \vert \nu)}.
$$
Here $I(\mu \vert \nu)$ is the Fisher information of $\mu$ with respect to $\nu$ defined in~\eqref{eq-fisher}.

If $\nu$ satisfies a $WJ(C)$ inequality then, in the notation ${\nabla\psi}\#\mu=\nu$,  
$$
W_2^2(\nu,\mu)\leq \frac1C J(\mu \vert \nu) \leq \frac1C \int (x-\n\psi)\cdot(\n\log{\mu}+A)d\mu
\leq
\frac1C W_2(\nu, \mu) \, \sqrt{I(\mu \vert \nu)}
$$
by the Cauchy-Schwarz inequality, as in Remark~\ref{rem:AGSOV}:
 
  \begin{eprop}\label{WJ-WI}
A $WJ$ inequality implies a $WI$ inequality with the same constant. 
 \end{eprop}
 
Let us now examine the link with the Talagrand~\eqref{WH} and logarithmic Sobolev~\eqref{LSI} inequalities.

 \begin{ecor}\label{cor:WJWHLSI}
1) A $WJ$ inequality implies a $WH$ inequality with the same constant. \\
 2) Assume that the probability measure $\nu=e^{-V}dx$ satisfies a $WJ(C)$ inequality and   $\nabla^2 V\geq \rho{\rm Id}_n$, for some $\rho\in\dR$. Then $\nu$ satisfies a logarithmic Sobolev with constant $C\left(1+\frac{\max(0, -\rho)}{2C}\right)^{-2}$.\\

 \end{ecor}

 \begin{eproof}
 1) By \cite[Th.~2.4]{glww}, a $WI$ inequality implies a $WH$ inequality with the same constant, so that the result comes from Proposition~\ref{WJ-WI}.\\

2) By  \cite[Th. 2]{ov00}, the following $HWI$ inequality holds: for all $\mu$
$$
H(\mu|\nu)\leq \sqrt{I(\mu|\nu)}W_2(\nu,\mu)+\frac{\rho_-}{2}W_2^2(\nu,\mu)= W_2(\nu, \mu) \, (\sqrt{I(\mu|\nu)}+\frac{\rho_-}{2} \, W_2(\nu,\mu) ).
$$
Here $\rho_- = \max(0, -\rho)$. As a $WJ(C)$ inequality implies both $WH(C)$ and $WI(C)$ inequalities, we get 
$$
H(\mu|\nu)\leq \sqrt{\frac{2}{C} H(\mu|\nu)} (1+\frac{\rho_-}{2C})\sqrt{I(\mu|\nu)}
$$
which ends the proof.
 \end{eproof}

Observe that in the uniformly convex case when $\rho >0$, then $\nu$ classically satisfies all WJ, WI, WH and logarithmic Sobolev inequalities with constant $C$. Moreover, under the assumption 2) with $C \geq \max( \rho, 0)$, then~\cite[Cor. 3.2]{ov00} ensures a log Sobolev inequality with constant $C(2-\rho/C)^{-1}$: for instance for $\rho  \leq 0$ it is worse than our constant (since then $C \geq \vert \rho \vert /2$).

Remark also, by~\cite{glwy} and~\cite{ov00}, that a $WI(C)$ or a $WH(C)$ inequality imply a Poincar\'e inequality with the same constant, hence providing an alternative proof to  Proposition~\ref{WJ-PI}. 

Observe finally that the general bound
$$
W_2(\mu_0, \nu) - W_2(\mu_t, \nu) \leq t^{1/2} H(\mu_0 \vert \nu)^{1/2}
$$
was obtained in~\cite[Remark 4.9]{cg06} for all $t$ and solutions $(\mu_t)_t$ to~\eqref{eq:FP} , hence directly proving that a uniform decay of the Wasserstein distance as in (\ref{eq-wasserstein-decay}) implies a $WH$ inequality with constant 
$\displaystyle \sup_{t >0} 2 \frac{(1-e^{-Ct})^2}{t} = 2 C \sup_{x >0} \frac{(1-e^{-x})^2}{x} \sim 0.8 \, C$ instead of $C$, which is optimal.

\medskip

We do not know whether a logarithmic Sobolev inequality, which implies a $WI$ inequality, also implies a $WJ$ inequality, or whether the converse holds without the curvature condition of Corollary~\ref{cor:WJWHLSI}.

 \subsection{Proof of Proposition~\ref{Amonotone}}
 
 We first state a general result on the map $A$:

\begin{elem}\label{lem:Vconvexe} 
Let $A$ be a $C^1$ monotone map on $\rr^n$ for which there exist two constants $R$ and $K>0$ such that $\nabla^S A(x) \geq K$ for all $\vert x \vert \geq R$. Then
$$(A (x) - A (y))  \cdot (x-y) \geq \frac{K}{3}\vert x-y \vert^2
$$
if $\vert x \vert \geq 2R$ or  $\vert y \vert \geq 2R$.
\end{elem}

\begin{eproof} 
Let $x$ and $y$ be fixed in $\rr^n$ with $\vert y \vert \geq 2 R$, and let us first write
\begin{eqnarray*}
(A(x)- A(y)) \cdot (x-y)
&=&
\int_0^1\nabla^S A (y+t(x-y))\,(x-y) \cdot (x-y) \, dt
\\
&=&
 r \int_0^r \nabla^S A(y+s \theta)\,\theta \cdot \theta \, ds
\end{eqnarray*}
for $x=y+r\theta$ with $r(=\vert x-y \vert)\geq 0$ and  $ \theta \in \mathbb S^{n-1}$.

{\bf 1.} If  $\{y+t(x-y);\,0\leq t \leq 1\}\cap \{z \in \rr^n;\,\vert z \vert \leq R\}=\emptyset$, then
$$
\int_0^1 \nabla^S A(y+t(x-y))\,(x-y) \cdot (x-y) \,dt \geq K \vert x-y \vert^2\geq
K \frac{\vert x-y\vert^2}{3}\cdot
$$

{\bf 2.} If  $\{y+t(x-y);\,0\leq t \leq 1\}\cap \{z \in \rr^n;\,\vert z \vert \leq R\}\not=\emptyset$, then let  $0\leq r_-\leq r_+$ such that
$$
\{y+s\theta;\,s \geq 0\}\cap \{z\in \rr^d;\,\vert z \vert \leq R\}=[ r_- \theta ,r_+\theta].
$$

Observe that
$$
 r_-=\vert y-(y+ r_\theta^- \theta)\vert\geq \inf\{\vert y-z\vert;\, \vert z \vert \leq 
 R\}=\vert y \vert-R
 $$
and
 $$
  r_+\leq \sup\{\vert y-z\vert;\, \vert z \vert \leq R\}=\vert y \vert+R
 $$
with $\vert y \vert \geq 2R$, so that
 $$
   r_- \geq \frac{  r_+}{3} \cdot
   $$

{\bf 2.1.} If  $r_-\leq r \leq   r_+$, then
$$
\int_0^r \nabla^S A(y+s \theta)\,\theta \cdot \theta\,ds \geq \int_0^{r_-} \nabla^S A(y+s \theta)\,\theta \cdot\theta\,ds\geq K r_-\geq
K \frac{r_+}{3}\geq
K \frac{r}{3}\cdot
$$

{\bf 2.2.} If  $r_+\leq r $, then
\begin{eqnarray*}
\int_0^r  \nabla^S A (y+s \theta)\,\theta \cdot \theta\,ds  \! \!
&\geq&
 \int_0^{r_-}  \nabla^S A(y+s \theta)\,\theta \cdot \theta\,ds + \int_{r_+}^r  \nabla^S A (y+s \theta)\,\theta \cdot \theta \,ds
\\
&\geq& \! \!
 K r_-+K(r-r_+)   \geq K \big(\frac{r_+}{3}+r-r_+ \big) 
=K \Big(\frac{r}{3} +\frac{2(r-r_+)}{3} \Big)\geq 
K \frac{r}{3}\cdot
\end{eqnarray*}

This concludes the argument, all cases being covered.
\end{eproof}

\bigskip

We now turn to the {\bf proof of Proposition~\ref{Amonotone}}. Let  $\varphi$  be a given strictly convex $C^2$ function on $\rr^n$. Let us recall that for the Hessian operator  
$$
\nabla^2 \varphi^* \, (\nabla\varphi(x))=(\nabla^2 \varphi(x))^{-1}
$$
and in particular
$$
\Delta\varphi^* \, (\nabla\varphi(x))=trace(\nabla^2\varphi(x))^{-1}
$$

Let  $X$ be    the subset of $\rr^n$ defined by
  $$
  X=\{x \in \rr^n,\, \vert x \vert \leq 2R, \, \vert \nabla\varphi (x) \vert \leq 2R\}.
  $$

 {\bf 1.} First of all, by monotonicity of $A$ and Lemma~\ref{lem:Vconvexe},
  \begin{multline*}
 \int_{\rr^n} (A(\nabla \varphi(x))- A(x)) \cdot (\nabla \varphi(x)-x) \,e^{-V(x)}\,dx 
 \\
\geq   \int_{\rr^n \setminus X} (A(\nabla \varphi(x))-A(x)) \cdot (\nabla \varphi(x)-x) \, e^{-V(x)}\,dx \geq \frac{K}{3}
  \int_{\rr^n\setminus X}\vert \nabla\varphi(x)-x\vert ^2 \, e^{-V(x)}\,dx.
 \end{multline*}
  
  \medskip
  
 {\bf 2.} On the other hand, for $\theta \in \dS^{n-1}$ we let  $R_{\theta}=\sup\{r\geq 0,\,r\theta \in X\}$. In particular $R_{\theta}\,\theta \in X$ and $R_{\theta} \leq 2R$.
 Then we let  $r_{\theta} \in [R_{\theta}, 3R]$ such that 
 $$
  \vert\nabla  \varphi(r_{\theta}\,\theta)-r_{\theta}\,\theta\vert =\inf\{  \vert \nabla\varphi(r\,\theta)-r\,\theta \vert , \,R_{\theta} \leq r \leq  3R\} .
    $$
    
In particular
$$
  \vert \nabla \varphi(r_{\theta}\,\theta)\vert \leq   \vert \nabla \varphi(r_{\theta}\,\theta)-r_{\theta}\,\theta \vert +\vert r_{\theta}\,\theta \vert  \leq   \vert \nabla \varphi(R_{\theta}\,\theta)-R_{\theta}\,\theta \vert +\vert r_{\theta}\,\theta \vert  \leq 2R+2R+3R = 7R
  $$ 
since   $\vert \nabla \varphi(R_{\theta}\,\theta)\vert  \leq 2R$ and $\vert R_{\theta}\,\theta \vert  \leq 2R$ for $R_{\theta}\,\theta \in X$.

Then, for $ r\, \theta \in X$ with $0\leq r \leq R_\theta \leq r_\theta$, let us write
$$
\nabla \varphi(r\theta)-r\theta=\nabla \varphi(r_\theta \theta)-r_\theta\,\theta+ \int_{r_\theta}^r [\nabla^2\varphi \, (s\theta)-I] \, \theta\,ds.
$$
Letting $H= \nabla^2 \varphi (s\theta)$ for notational convenience,  we decompose as
$$
[H-I] \ \theta = [H^{\frac{1}{2}}-H^{-\frac{1}{2}}]H^\frac{1}{2} \, \theta
$$
so that
\begin{eqnarray*}
\Big\vert \int_r^{r_\theta} [H-I] \theta \,dt  \Big\vert^2
&\leq&
\Big( \int_r^{r_\theta}\vert H^\frac{1}{2}-H^{-\frac{1}{2}}\vert \vert H^{\frac{1}{2}} \theta \vert \,ds \Big)^2
\\
&\leq& \int_r^{r_\theta}\vert H^\frac{1}{2}-H^{-\frac{1}{2}}\vert ^2\,e^{-V(s\theta)}\,ds \; \int_r^{r_\theta }\vert H^{\frac{1}{2}} \theta\vert^2 e^{+V(s\theta)}\,\,ds.
\end{eqnarray*}
by the H\"older inequality. But
\begin{eqnarray*}
\vert H^\frac{1}{2}-H^{-\frac{1}{2}}\vert ^2
&=&
\sup_x \frac{\vert [H^\frac{1}{2}-H^{-\frac{1}{2}}]x\vert ^2}{\vert x \vert^2}
=\sup_x \frac{([H-2I+H^{-1}]x) \cdot x}{\vert x \vert^2}
\\
&\leq&
  trace(H-2I+H^{-1}) =\Delta\varphi (s\theta) -2n+(\Delta\varphi^*)(\nabla \varphi ((s\theta)).
\end{eqnarray*}
since the eigenvalues of $H-2I+H^{-1}$ are non-negative. Moreover
$$
 \vert H^{\frac{1}{2}} \, \theta \vert^2
=
 (H^{\frac{1}{2}} \, \theta) \cdot (H^{\frac{1}{2}} \, \theta )=H \, \theta \cdot \theta.
 $$
 Hence
\begin{multline*}
\vert \nabla \varphi(r\,\theta)-r\, \theta\vert^2 \leq 2 \, \vert \nabla \varphi(r_\theta\,\theta)-r_\theta\,\theta\vert^2
\\
+2 \int_r^{r_\theta}(\Delta\varphi(s\,\theta)-2n+\Delta\varphi^*(\nabla \varphi(s\,\theta))\,e^{-V(s\,\theta)}\,ds \;
\,\int_r^{r_\theta} (H \, \theta) \cdot \theta \, e^{+V(s\theta)}\,ds 
\end{multline*}
where
 $$
 \int_r^{r_\theta} ( H \, \theta) \cdot \theta\,ds 
 =(\nabla \varphi(r_\theta\,\theta)-\nabla \varphi(r\,\theta)) \cdot \theta
 \leq 
 \vert \nabla \varphi(r_\theta\theta)\vert +\vert \nabla \varphi (r\,\theta)\vert \leq 9R
 $$
for  $r\,\theta \in X$. Hence
\begin{eqnarray*}
\int_{X, \vert x \vert \leq 2R}\vert \nabla \varphi(x)-x\vert^2e^{-V(x)}\,dx 
&=&
\int_{\dS^{n-1}}\int_0^{R_\theta} r^{n-1}\vert \nabla \varphi(r\,\theta)-r\theta\vert^2\,e^{-V(r \theta)}\,dr \,d\theta
\\
&\leq&
 2 \int_{\dS^{n-1}}\int_0^{R_\theta} r^{n-1}\vert \nabla \varphi(r_\theta\,\theta)-r_\theta\theta\vert^2\,e^{-V(r \theta)}\,dr \,d\theta
\end{eqnarray*}
$$
+18R\,e^{\sup\{V(x);\,\vert x \vert \leq 2R\}} \,\int_{\dS^{n-1}}\int_0^{R_\theta} r^{n-1}
\int_r^{r_\theta}(\Delta\varphi(s\theta)-2n+\Delta\varphi^*(\nabla \varphi (s\theta))\,e^{-V(s\theta)}\,ds \,dr\,d\theta.
$$

But
\begin{eqnarray*}
&&
\int_{\dS^{n-1}}\int_0^{R_\theta}r^{n-1}
\int_r^{r_\theta}(\Delta\varphi(s\theta)-2n+\Delta\varphi^*(\nabla \varphi (s\theta))\,e^{-V(s\theta)}\,ds \,dr\,d\theta
\\
&&
\leq 
\int_{\dS^{n-1}}\int_0^{R_\theta}
\int_r^{r_\theta}s^{n-1}(\Delta\varphi(s\theta)-2n+\Delta\varphi^*(\nabla \varphi (s\theta))\,e^{-V(s\theta)}\,ds \,dr\,d\theta
\\
&&
\leq 2R
\int_{\dS^{n-1}}
\int_0^{3R}s^{n-1}(\Delta\varphi(s\theta)-2n+\Delta\varphi^*(\nabla \varphi (s\theta))\,e^{-V(s\theta)}\,ds \,d\theta
\\
&&
 =2R
\int_{\vert x \vert \leq 3R}
(\Delta\varphi(x)-2n+\Delta\varphi^*(\nabla \varphi (x)))\,e^{-V(x)}\,dx.
\end{eqnarray*}
 
Hence
 \begin{eqnarray}
&&
\int_{X, \vert x \vert \leq 2R}\vert \nabla \varphi(x)-x\vert^2e^{-V(x)}\,dx 
\nonumber\\
&\leq&
 2 \, e^{-\inf\{V(x);\,\vert x \vert \leq 2R\}} \frac{(2R)^{n}}{n}\,\, \int_{\dS^{n-1}}\vert \nabla \varphi(r_\theta\,\theta)-r_\theta\theta\vert^2 \,d\theta
\nonumber\\
&&
+ \, 18R\,e^{\sup\{V(x);\,\vert x \vert \leq 2R\}} \,2R
\int_{\vert x \vert \leq 3R}(\Delta\varphi(x)-2n+\Delta\varphi^*(\nabla \varphi (x)))\,e^{-V(x)}\,dx.
 \label{Acvx1}
  \end{eqnarray}

Moreover, by Lemma \ref{lem:Vconvexe}  and the definition of $r_\theta$,
\begin{eqnarray*}
&& \int_{\vert x \vert \leq 3R} (A(\nabla \varphi(x))- A(x)) \cdot (\nabla \varphi(x)-x) \, e^{-V(x)}\,dx
\\
 && \geq
  \int_{2R \leq\vert x \vert \leq 3R} (A(\nabla \varphi(x))- A(x)) \cdot (\nabla \varphi(x)-x) \, e^{-V(x)}\,dx
\\
&&
 \geq 
 \frac{K}{3}\int_{2R \leq\vert x \vert \leq 3R}\vert \nabla \varphi(x)-x\vert^2 e^{-V(x)}\,dx
 \geq 
 \frac{K}{3}e^{-\sup\{V(x);\,\vert x \vert \leq 3R\}}\int_{2R \leq\vert x \vert \leq 3R}\vert \nabla \varphi(x)-x\vert^2 \,dx
 \\
 &&
 = \frac{K}{3}e^{-\sup\{V(x);\,\vert x \vert \leq 3R\}}\int_{2R} ^{3R}r^{n-1}\int_{\dS^{n-1}}\vert \nabla \varphi(r\,\theta)-r\,\theta\vert^2 \,dr\,d\theta
 \\
 &&
 \geq 
 \frac{K}{3}e^{-\sup\{V(x);\,\vert x \vert \leq 3R\}}\int_{2R}^{3R}r^{n-1}\int_{\dS^{n-1}}\vert \nabla \varphi(r_\theta\,\theta)-r_\theta\,\theta\vert^2 \,dr\,d\theta.
\end{eqnarray*}
 
 Hence there exists a constant $C$ such that
 \begin{equation}\label{Acvx2}
C \int_{\dS^{n-1}}\vert\nabla \varphi(r_\theta \theta)-r_\theta\,\theta\vert^2\, \, d\theta
 \leq  \int_{\vert x \vert \leq 3R} (A(\nabla \varphi(x))- A(x)) \cdot (\nabla \varphi(x)-x) \, e^{-V(x)}\,dx.
\end{equation}

It follows from~\eqref{Acvx1} and~\eqref{Acvx2}  that
\begin{eqnarray*}
C  \int_{X, \vert x \vert \leq 2R}\vert \nabla  \varphi(x)-x\vert^2\, e^{-V(x)}\,dx 
  &\leq &
  \int_{\vert x \vert \leq 3R} (A(\nabla \varphi(x)- A (x)) \cdot (\nabla \varphi(x)-x) \, e^{-V(x)}\,dx  
  \\
&&  +\int_{\vert x \vert \leq 3R}(\Delta\varphi(x)-2n+\Delta\varphi^*(\nabla \varphi(x))\,e^{-V(x)}\,dx.
\end{eqnarray*}

Moreover 
$$
 \int_{\vert x \vert \leq 3R} (A(\nabla \varphi(x))- A(x)) \cdot (\nabla \varphi(x)-x) \, e^{-V(x)}\,dx
 \geq 
 \frac{K}{3}\int_{2R \leq\vert x \vert \leq 3R}\vert \nabla \varphi(x)-x\vert^2 e^{-V(x)}\,dx, 
$$
so that
\begin{eqnarray*}
C  \int_{X}\vert \nabla  \varphi(x)-x\vert^2\, e^{-V(x)}\,dx 
  &\leq &
\int_{\vert x \vert \leq 3R} (A(\nabla \varphi(x)- A(x)) \cdot (\nabla \varphi(x)-x) \, e^{-V(x)}\,dx  
  \\
&&  +\int_{\vert x \vert \leq 3R}(\Delta\varphi(x)-2n+\Delta\varphi^*(\nabla \varphi(x))\,e^{-V(x)}\,dx.
\end{eqnarray*}   
   
Finally the last two integrands are non-negative maps, so we can bound from above these last two integrals on  the set  $\{\vert x \vert \leq 3R\}$ by the corresponding integrals on the whole~$\rr^n$.

\medskip

{\bf 3.} We conclude the proof of Proposition~\ref{Amonotone} by adding the estimates in {\bf 1.} and {\bf 2.} 

\noindent

\bigskip

{\bf Acknowledgements.} The authors are grateful to G. Carlier and A. Figalli for enlighting discussion. They thank a 
referee for a careful reading of the manuscript and most relevant comments, references and questions which helped improve the presentation of the paper. This research was supported in part by the ANR project EVOL.

\newcommand{\etalchar}[1]{$^{#1}$}


\end{document}